\documentclass[english]{amsart}
 
\usepackage{amsmath}

\usepackage{amssymb}

\usepackage{amscd}
\usepackage[all,cmtip,line]{xy} \usepackage{enumerate}

\newcommand{\Wconj}{W^?}
\newcommand{\Lcrys}{L_{\operatorname{crys}}}
\newcommand{\Lflat}{L_{\operatorname{flat}}}

\newcommand{\ur}{\mathrm{ur}}

\newcommand{\into}{\hookrightarrow}

\newcommand{\dd}{\operatorname{dd}}

\newcommand{\BrMod}{\operatorname{BrMod}}

\newlength{\ownl}

\newcommand{\Art}{{\operatorname{Art}}}

\newcommand{\Ext}{{\operatorname{Ext}}}
\newcommand{\Fil}{{\operatorname{Fil}}}

\newcommand{\Frob}{{\operatorname{Frob}}}

\newcommand{\Gal}{{\operatorname{Gal}}}

\newcommand{\Hom}{{\operatorname{Hom}}}

\newcommand{\GL}{\operatorname{GL}}

\newcommand{\ab}{{\operatorname{ab}}}

\newcommand{\st}{{\operatorname{st}}}

\newcommand{\F}{{\mathbb{F}}}

\newcommand{\M}{{\mathcal{M}}}
\newcommand{\N}{{\mathcal{N}}}

\newcommand{\Q}{{\mathbb{Q}}}

\newcommand{\Z}{{\mathbb{Z}}}

\newcommand{\cG}{\mathcal{G}}

\newcommand{\cM}{\mathcal{M}}
\newcommand{\cN}{\mathcal{N}}
\newcommand{\cO}{\mathcal{O}}
\newcommand{\cP}{\mathcal{P}}

\newcommand{\barK}{\overline{{K}}}

\newcommand{\barQQ}{\overline{{\Q}}}

\newcommand{\tv}{{\widetilde{{v}}}}

\newcommand{\epsilonbar    }{\overline{\epsilon}}

 \newcommand{\thetabar    }{\overline{\theta}}

 \newcommand{\rhobar   }{{\overline{\rho}}}

 \newcommand{\chibar   }{\overline{\chi}}

 \newcommand{\omegat   }{\widetilde{\omega}}

\def\RCS$#1: #2 ${\expandafter\def\csname RCS#1\endcsname{#2}}
\RCS$Revision: 60 $
\RCS$Date: 2011-06-24 13:12:16 +0100 (Fri, 24 Jun 2011) $

\newcommand{\bigO}{\mathcal{O}}

\newcommand{\bb}{\mathbb}
\newcommand{\mc}{\mathcal}
\newcommand{\mf}{\mathfrak}

\newcommand{\rbar}{\bar{r}}

\newcommand{\HT}{\operatorname{HT}}

 \newcommand{\Qp}{\Q_p}
 
\newcommand{\Qpbar}{{\overline{\Q}_p}}
\newcommand{\Zpbar}{\overline{\Z}_p}

\newcommand{\Fpbar}{{\overline{\F}_p}}

\newcommand{\Fp}{{\F_p}}

\newcommand{\Sym}{\operatorname{Sym}}

\usepackage{hyperref} 
\usepackage{amsthm}

\newtheorem{ithm}{Theorem}
\newtheorem{thm}{Theorem}[subsection]

\newtheorem{cor}[thm]{Corollary}
 
\newtheorem{lem}[thm]{Lemma}
\newtheorem{prop}[thm]{Proposition}
\newtheorem{conj}[thm]{Conjecture} \theoremstyle{definition}
 \theoremstyle{definition}
\newtheorem{defn}[thm]{Definition} \theoremstyle{remark}
\newtheorem{rem}[thm]{Remark}
\newtheorem{remark}[thm]{Remark}

\numberwithin{equation}{subsection}

\theoremstyle{definition}

\begin{document}
\title[]{Crystalline extensions and the weight part of Serre's conjecture}

\author{Toby Gee} \email{gee@math.northwestern.edu} \address{Department of
  Mathematics, Northwestern University}\author{Tong
  Liu}\email{tongliu@math.purdue.edu}\address{Department of
  Mathematics, Purdue University}\author{David Savitt} \email{savitt@math.arizona.edu}
\address{Department of Mathematics, University of Arizona}
\thanks{The authors were partially
  supported by NSF grants DMS-0841491, DMS-0901360, and DMS-0901049 respectively.}  
\begin{abstract}Let $p>2$ be prime. We complete the proof of the
  weight part of Serre's conjecture for rank two unitary groups for
  mod $p$ representations in the totally ramified case, by proving
  that any weight which occurs is a predicted weight. Our methods are
  a mixture of local and global techniques, and in the course of the
  proof we establish some purely local results on crystalline
  extension classes.
\end{abstract}
\maketitle
\tableofcontents

\section{Introduction}\label{sec: introduction} The weight
part of generalisations of Serre's conjecture has seen significant
progress in recent years, particularly for (forms of) $\GL_2$. Conjectural
descriptions of the set of Serre weights were made in increasing
generality by \cite{bdj}, \cite{MR2430440} and \cite{GHS}, and cases
of these conjectures were proved in \cite{geebdj} and
\cite{geesavitttotallyramified}. Most recently, significant progress
was made towards completely establishing the conjecture for rank two
unitary groups in \cite{blggU2}. We briefly recall this result. Let
$p>2$ be prime, let $F$ be a CM field, and let
$\rbar:G_F\to\GL_2(\Fpbar)$ be a modular representation (see
\cite{blggU2} for the precise definition of ``modular'', which is in
terms of automorphic forms on compact unitary groups). There is a
conjectural set $\Wconj(\rbar)$ of Serre weights in which $\rbar$ is
predicted to be modular, which is defined in Section \ref{sec: serre
  weight definitions} below, following \cite{GHS}. Then the main
result of \cite{blggU2} is that under mild technical hypotheses,
$\rbar$ is modular of every weight in $\Wconj(\rbar)$.

It remains to show that if $\rbar$ is modular of some weight, then
this weight is contained in $\Wconj(\rbar)$. It had been previously
supposed that this was the easier direction; indeed, just as in the
classical case, the results of
\cite{blggU2} reduce the weight part of Serre's conjecture for these
unitary groups to a purely local problem in $p$-adic Hodge
theory. However, this problem has proved to be difficult,
and so far only fragmentary results are
known. In the present paper we resolve the problem in the totally
ramified case, so that in combination with \cite{blggU2} we resolve
the weight part of Serre's conjecture in this case, proving the
following Theorem (see Theorem \ref{thm: the main result, modular if
  and only if predicted}).
\begin{ithm}
  \label{thm: intro: the main result, modular if and only if predicted}Let
  $F$ be an imaginary CM field with maximal totally real subfield
 ~$F^+$, and suppose that $F/F^+$ is unramified at all finite places,
  that every place of $F^+$ dividing $p$ splits completely in $F$,
  that $\zeta_p\notin F$, and that $[F^+:\Q]$ is even. Suppose that
  $p>2$, and that $\rbar:G_F\to\GL_2(\Fpbar)$ is an irreducible
  modular representation with split ramification such that
  $\rbar(G_{F(\zeta_p)})$ is adequate. Assume that for each place $w|p$
  of $F$, $F_w/\Qp$ is totally ramified.

 Let $a\in(\Z^2_+)_0^S$ be a Serre weight. Then
 $a_w\in\Wconj(\rbar|_{G_{F_w}})$ if and only if $\rbar$ is modular of
 weight $a$.
\end{ithm}(See the body of the paper, especially Section~\ref{ss:global}, for any unfamiliar notation and
terminology.) While \cite{blggU2} reduced this result to a purely
local problem, our methods are not purely local; in fact we use the
main result of \cite{blggU2}, together with potential automorphy
theorems, as part of our proof.

In the case that $\rbar|_{G_{F_w}}$ is semisimple for each place
$w|p$, the result was established (in a slightly different setting) in
\cite{geesavitttotallyramified}. The method of proof was in part
global, making use of certain potentially Barsotti-Tate lifts to
obtain conditions on $\rbar|_{G_{F_w}}$. We extend this analysis in
the present paper to the case that $\rbar|_{G_{F_w}}$ is reducible but
non-split,
obtaining conditions on the extension classes that can occur; we show
that (other than in one exceptional case) they lie in a certain set $\Lflat$, defined in terms of finite
flat models.

In the case that $\rbar|_{G_{F_w}}$ is reducible the definition of
$\Wconj$ also depends on the extension class; it is required to lie in
a set $\Lcrys$, defined in terms of reducible crystalline lifts with
specified Hodge-Tate weights. To complete the proof, one must show
that $\Lcrys=\Lflat$. An analogous result was proved in generic
unramified cases in section 3.4 of \cite{geebdj} by means of explicit
calculations with Breuil modules; our approach here is less direct,
but has the advantage of working in non-generic cases, and requires
far less calculation.

We use a global argument to show that
$\Lcrys\subset\Lflat$. Given a class in $\Lcrys$, we use potential
automorphy theorems to realise the corresponding local representation
as part of a global modular representation, and then apply the main
result of \cite{blggU2} to show that this representation is modular of
the expected weight. Standard congruences between automorphic forms
then show that this class is also contained in $\Lflat$.

To prove the converse inclusion, we make a study of different finite
flat models to show that $\Lflat$ is contained in a vector space of
some dimension $d$. A standard calculation shows that $\Lcrys$
contains a space of dimension $d$, so equality follows. As a
byproduct, we show that both $\Lflat$ and $\Lcrys$ are vector
spaces. We also show that various spaces defined in terms of
crystalline lifts are independent of the choice of lift (see Corollary
\ref{cor: independence of lift for H^1_f}). The analogous property was
conjectured in the unramified case in \cite{bdj}.

It is natural to ask whether our methods could be extended to handle
the general case, where $F_w/\Qp$ is an arbitrary
extension. Unfortunately, this does not seem to be the case, because
in general the connection between being modular of some weight and
having a potentially Barsotti-Tate lift of some type is less direct. We expect that our methods could be used to reprove the results of
section 3.4 of \cite{geebdj}, but we do not see how to extend them to
cover the unramified case completely.

We now explain the structure of the paper. In Section \ref{sec: serre
  weight definitions} we recall the definition of~$\Wconj$, and the
global results from \cite{blggU2} that we will need. In Section \ref{sec:local to
  global} we recall a potential automorphy result from \cite{geekisin}, allowing us to
realise a local mod $p$ representation globally. Section \ref{sec:
  congruences to weight 0} contains the definitions of the spaces
$\Lcrys$ and $\Lflat$ and the proof that $\Lcrys\subset\Lflat$, and in
Section \ref{sec: finite flat
  models} we carry out the necessary calculations with Breuil modules
to prove our main local results. Finally, in section \ref{sec: global
  consequences} we combine our local results with the techniques of
\cite{geesavitttotallyramified} and the main result of \cite{blggU2}
to prove Theorem \ref{thm: intro: the main result, modular if and only if predicted}.

\subsection{Notation}If $M$ is a field, we let $G_M$ denote its
absolute Galois group. Let~$\epsilon$ denote the $p$-adic cyclotomic
character, and $\bar{\epsilon}$ the mod $p$ cyclotomic character.
If~$M$ is a global field and $v$ is a place of $M$, let $M_v$ denote
the completion of $M$ at $v$.   If
~$M$ is a finite extension of $\bb{Q}_l$ for some $l$, we write $I_M$
for the inertia subgroup of~$G_M$. If $R$ is a local ring we write
$\mf{m}_{R}$ for the maximal ideal of $R$.

 Let $K$ be a finite extension of $\Qp$, with ring of integers $\cO_K$
 and residue field~$k$.   We write $\Art_K:K^\times\to W_K^{\ab}$ for
the isomorphism of local class field theory, normalised so that
uniformisers correspond to geometric Frobenius elements.  For each $\sigma\in \Hom(k,\Fpbar)$ we
define the fundamental character $\omega_{\sigma}$ corresponding
to~$\sigma$ to be the composite $$\xymatrix{I_K \ar[r] & W_K^{\ab} \ar[r]^{\Art_K^{-1}} &
  \bigO_{K}^{\times}\ar[r] & k^{\times}\ar[r]^{\sigma} &
  \Fpbar^{\times}.}$$
In the case that $k\cong\Fp$, we will sometimes write $\omega$ for
$\omega_\sigma$. Note that in this case we have $\omega^{[K:\Qp]}=\epsilonbar$.

We fix an algebraic closure $\barK$ of $K$.  If $W$ is a de Rham representation of $G_K$ over
$\barQQ_p$ and $\tau$ is an embedding $K \into \barQQ_p$ then the multiset
$\HT_\tau(W)$ of Hodge-Tate weights of $W$ with respect to $\tau$ is
defined to contain the integer $i$ with multiplicity $$\dim_{\barQQ_p} (W
\otimes_{\tau,K} \widehat{\barK}(-i))^{G_K},$$ with the usual notation
for Tate twists.  Thus for example
$\HT_\tau(\epsilon)=\{ 1\}$.

\section{Serre weight conjectures: definitions}\label{sec: serre
  weight definitions}\subsection{Local definitions}We begin by recalling some
generalisations of the weight part of Serre's conjecture. We begin
with some purely local definitions. Let $K$ be a finite totally ramified extension of
$\Qp$ with absolute ramification index $e$, and let $\rhobar:G_K\to\GL_2(\Fpbar)$ be a continuous
representation.
\begin{defn}
  A \emph{Serre weight} is an irreducible $\Fpbar$-representation of
  $\GL_2(\Fp)$. Up to isomorphism, any such representation is of the
  form \[F_a:=\det{}^{a_{2}}\otimes\Sym^{a_{1}-a_{2}}\F_p^2\otimes_{\F_p}\Fpbar\]
  where $0\le a_{1}-a_{2}\le p-1$.  We also use the term Serre weight
  to refer to the pair $a = (a_1,a_2)$.
\end{defn}
We say that two Serre weights $a$ and $b$ are \emph{equivalent} if and only if
$F_a\cong F_b$ as representations of $\GL_2(\Fp)$. This is equivalent
to demanding that we
have $a_{1}-a_{2}=b_{1}-b_{2}$ and $a_2\equiv b_2\pmod{p-1}$.

We write $\Z^2_+$ for the set of pairs of integers $(n_1,n_2)$ with
$n_1\ge n_2$, so that a Serre weight $a$ is by definition an element
of $\Z^2_+$. We say that an element
$\lambda\in(\Z^2_+)^{\Hom_{\Qp}(K,\Qpbar)}$ is a \emph{lift} of a weight
$a\in\Z^2_+$ if there is an element $\tau\in\Hom_{\Qp}(K,\Qpbar)$ such that
$\lambda_{\tau}=a$, and for all other $\tau'\in\Hom_{\Qp}(K,\Qpbar)$ we have
$\lambda_{\tau'}=(0,0)$.

\begin{defn}
  \label{defn: Galois representation of Hodge type some weight}Let
  $K/\Qp$ be a finite extension, let
  $\lambda\in(\Z^2_+)^{\Hom_{\Qp}(K,\Qpbar)}$, and let
  $\rho:G_K\to\GL_2(\Qpbar)$ be a de Rham representation. Then we say
  that $\rho$ has \emph{Hodge type} $\lambda$ if for each
  $\tau\in\Hom_{\Qp}(K,\Qpbar)$ we have $\HT_\tau(\rho)=\{\lambda_{\tau,1}+1,\lambda_{\tau,2}\}$.
\end{defn}

Following \cite{GHS} (which in turn follows \cite{bdj} and \cite{MR2430440}), we define an explicit
set of Serre weights $\Wconj(\rhobar)$.
\begin{defn}
  \label{defn: W? niveau 1}If $\rhobar$ is reducible, then a Serre
  weight $a\in\Z^2_+$ is in $\Wconj(\rhobar)$ if
  and only if $\rhobar$ has a crystalline lift of the
  form \[ \begin{pmatrix}\psi_1&*\\ 0& \psi_2
  \end{pmatrix}\] which has Hodge type $\lambda$ for some lift
  $\lambda\in(\Z^2_+)^{\Hom_{\Qp}(K,\Qpbar)}$ of $a$. In particular, if $a\in \Wconj(\rhobar)$ then by  Lemma 6.2 of \cite{geesavitttotallyramified} it is necessarily the case that  there is a decomposition
  $\Hom(\Fp,\Fpbar)=J\coprod J^c$ and  an integer
  $0\le \delta\le e-1$ such that \[\rhobar|_{I_K}\cong
  \begin{pmatrix} \omega^{\delta}
 \prod_{ \sigma\in
      J}\omega_{\sigma}^{a_{1}+1}\prod_{\sigma\in
      J^c}\omega_\sigma^{a_{2}}&*\\ 0& \omega^{e-1-\delta} \prod_{\sigma\in
      J^c}\omega_\sigma^{a_{1}+1}\prod_{\sigma\in
      J}\omega_\sigma^{a_{2}}.   \end{pmatrix}\]
\end{defn}
We remark that while it may seem strange to consider the single
element set $\Hom(\Fp,\Fpbar)$, this notation will be convenient for us.

\begin{defn}
  \label{defn: W? niveau 2}
 Let $K'$ denote the quadratic unramified
extension of $K$ inside~$\barK$, with residue field
$k'$ of order $p^2$.

 If $\rhobar$ is irreducible, then a Serre
  weight $a\in\Z^2_+$ is in $\Wconj(\rhobar)$ if
  and only if there is a subset $J\subset\Hom(k',\Fpbar)$ of size $1$,
  and  an integer $0\le
  \delta\le e-1$ such that if we write
  $\Hom(k',\Fpbar)=J\coprod J^c$, then  \[\rhobar|_{I_K}\cong
  \begin{pmatrix}\prod_{\sigma\in
      J}\omega_{\sigma}^{a_{1}+1+\delta}\prod_{\sigma\in
      J^c}\omega_\sigma^{a_{2}+e-1-\delta}&0\\ 0& \prod_{\sigma\in
      J^c}\omega_\sigma^{a_{1}+1+\delta}\prod_{\sigma\in
      J}\omega_\sigma^{a_{2}+e-1-\delta}
  \end{pmatrix}.\]
\end{defn}
We remark that by Lemma 4.1.19 of \cite{blggU2}, if
$a\in\Wconj(\rhobar)$ and $\rhobar$ is irreducible then
$\rhobar$ necessarily has a crystalline lift of Hodge type $\lambda$ for any lift
  $\lambda\in(\Z^2_+)^{\Hom_{\Qp}(K,\Qpbar)}$ of $a$. Note also that if $a$
  and $b$ are equivalent and $a\in\Wconj(\rhobar)$ then $b\in\Wconj(\rhobar)$.
  \begin{remark}\label{rem: conjectured weights independent of
      unramified twist}
    Note that if $\thetabar: G_K\to\Fpbar^\times$ is an unramified character, then
    $\Wconj(\rbar)=\Wconj(\rbar\otimes\thetabar)$. \end{remark}
  \subsection{Global conjectures}\label{ss:global} The point of the local definitions
  above is to allow us to formulate global Serre weight
  conjectures. Following \cite{blggU2}, we work with rank two unitary
  groups which are compact at infinity. As we will not need to make
  any arguments that depend on the particular definitions made in
  \cite{blggU2}, and our main results are purely local, we simply
  recall some notation and basic properties of the definitions,
  referring the reader to \cite{blggU2} for precise formulations.

We emphasise that our conventions for Hodge-Tate weights are the
opposite of those of \cite{blggU2}; for this reason, we must introduce
a dual into the definitions.

Fix an imaginary CM field $F$, and let $F^+$ be its maximal totally
real subfield. We assume that each prime of $F^+$ over $p$ has residue
field $\F_p$ and splits in $F$. We define a global notion of Serre weight by taking a
product of local weights in the following way.

\begin{defn}
  \label{defn:global-serre-wts}
Let $S$ denote the set
of places of $F$ above $p$.  If $w \in S$ lies over a place $v$ of
$F^+$, write $v = w w^c$.  Let
$(\Z^2_+)_0^{S}$ denote the subset of
$(\Z^2_+)^{S}$ consisting of elements $a = (a_w)_{w \in S}$ such
that $a_{w,1}+a_{w^c,2}=0$ for all $w\in S$. We say that an
element  $a\in(\Z^2_+)_0^{S}$ is a \emph{Serre
 weight} if for each $w|p$ we
have \[p-1\ge a_{w,1}-a_{w,2}.\]
\end{defn}

Let $\rbar:G_F\to\GL_2(\Fpbar)$ be a continuous irreducible
representation. Definition 2.1.9 of \cite{blggU2} states what it
means for $\rbar$ to be modular, and more precisely for $\rbar$ to be
modular of some Serre weight $a$; roughly speaking, $\rbar$ is modular
of weight $a$ if there is a cohomology class on some unitary group
with coefficients in the local system corresponding to $a$ whose
Hecke eigenvalues are determined by the characteristic polynomials of
$\rbar$ at Frobenius elements. Since our conventions for Hodge-Tate
weights are the opposite of those of \cite{blggU2}, we make the
following definition.

\begin{defn}
   Suppose that $\rbar:G_F\to\GL_2(\Fpbar)$ is a continuous
  irreducible modular representation. Then we say that $\rbar$ \emph{is modular
  of weight} $a\in(\Z^2_+)_0^S$ if
$\rbar^\vee$ is modular of weight $a$ in the sense of Definition 2.1.9
of \cite{blggU2}.
\end{defn} We  globalise the definition of the set
$\Wconj(\rhobar)$ in the following natural fashion.
\begin{defn}
  If $\rbar:G_F\to\GL_2(\Fpbar)$ is a continuous representation, then
  we define $\Wconj(\rbar)$ to be the set of Serre weights
  $a\in(\Z^2_+)_0^S$ such that for each
  place $w|p$ the corresponding Serre weight
  $a_w\in\Z^2_+$ is an element of
  $\Wconj(\rbar|_{G_{F_w}})$.
\end{defn}

One then has the following conjecture.
\begin{conj}\label{conj: global Serre weight explicit conjecture}
  Suppose that $\rbar:G_F\to\GL_2(\Fpbar)$ is a continuous irreducible
  modular representation, and that
  $a\in(\Z^2_+)_0^S$ is a Serre
  weight. Then $\rbar$ is modular of weight $a$ if and only if
  $a\in\Wconj(\rbar)$.
\end{conj}
If $\rbar:G_F\to\GL_2(\Fpbar)$ is a continuous representation,
  then we say that $\rbar$ has \emph{split ramification} if any finite
  place of $F$ at which $\rbar$ is ramified is split over $F^+$.   The
  following result is
Theorem 5.1.3 of \cite{blggU2}, one of the
main theorems of that paper, in the special case where $F_w/\Qp$ is
totally ramified for all $w|p$.  (Note that in \cite{blggU2}, the set
of weights $\Wconj(\rbar)$ is referred to as
$W^{\operatorname{explicit}}(\rbar)$.)
\begin{thm}
  \label{thm: explicit local lifts implies Serre
    weight}Let $F$ be an imaginary CM field with maximal totally real subfield~$F^+$. Assume that $\zeta_p\notin F$, that $F/F^+$ is unramified at all finite places,
  that every place of $F^+$ dividing $p$ has residue field $\F_p$ and splits completely in $F$,
  and that $[F^+:\Q]$ is even. Suppose that $p>2$, and that
  $\rbar:G_F\to\GL_2(\Fpbar)$ is an irreducible modular
  representation with split ramification. Assume that $\rbar(G_{F(\zeta_p)})$ is adequate.

 Let $a\in(\Z^2_+)_0^S$ be a
  Serre weight. Assume that $a\in \Wconj(\rbar)$. Then $\rbar$ is
  modular of weight $a$.
\end{thm}
Here \emph{adequacy} is a group-theoretic condition, introduced in
\cite{jack}, that for subgroups of $\GL_2(\Fpbar)$ with $p > 5$ is
equivalent to the usual condition that $\rbar|_{G_{F(\zeta_p)}}$ is irreducible.  For a precise
definition we refer the reader to Definition A.1.1 of \cite{blggU2}.
We  also remark that the hypotheses that  $F/F^+$ is unramified at all finite places,
  that every place of $F^+$ dividing $p$ splits completely in $F$,
  and that $[F^+:\Q]$ is even, are in fact part of the definition of
  ``modular'' made in \cite{blggU2}.)
  Theorem~\ref{thm: explicit local lifts implies Serre
    weight}  establishes one direction of Conjecture \ref{conj: global Serre
    weight explicit conjecture}, and we are left with the problem of
  ``elimination,'' i.e., the problem of proving that if $\rbar$ is
  modular of weight $a$, then $a\in \Wconj(\rbar)$.
  We believe that this problem should have a purely local resolution,
  as we now explain.

The key point is the relationship between being
  modular of weight $a$, and the existence of certain de Rham lifts of
  the local Galois representations $\rbar|_{G_{F_w}}$, $w|p$. The link
  between these properties is provided by local-global compatibility
  for the Galois representations associated to the automorphic
  representations under consideration; rather than give a detailed
  development of this connection, for which see \cite{blggU2}, we
  simply summarise the key results from \cite{blggU2} that we will
  use. The
  following is Corollary 4.1.8 of \cite{blggU2}.

\begin{prop}
  \label{prop: modular of some weight implies crystalline lifts
    exist}Let $F$ be an imaginary CM field with maximal totally real
  subfield $F^+$, and suppose that $F/F^+$ is unramified at all finite
  places, that every place of $F^+$ dividing $p$ has residue field
  $\F_p$ and splits completely in
  $F$, and that $[F^+:\Q]$ is even. Suppose that $p>2$, and that
  $\rbar:G_F\to\GL_2(\Fpbar)$ is an irreducible modular representation
  with split ramification. Let
  $a\in(\Z^2_+)_0^S$ be a Serre
  weight. If $\rbar$ is modular of weight $a$, then for each place
  $w|p$ of $F$, there is a crystalline representation
  $\rho_w:G_{F_w}\to\GL_2(\Qpbar)$ lifting $\rbar|_{G_{F_w}}$, such
  that $\rho_w$ has Hodge type $\lambda_w$ for some lift
  $\lambda_w\in(\Z^2_+)^{\Hom_{\Qp}(F_w,\Qpbar)}$ of $a$.
\end{prop}

We stress that Proposition~\ref{prop: modular of some weight implies crystalline lifts
    exist} does not already complete the proof of Conjecture \ref{conj: global Serre
    weight explicit conjecture}, because the representation $\rho_w$
  may be irreducible (compare with Definition~\ref{defn: W? niveau 1}).
However, in light of this result, it is natural to make the following
purely local conjecture, which together with Theorem \ref{thm:
  explicit local lifts implies Serre weight} would essentially resolve
Conjecture \ref{conj: global Serre weight explicit conjecture}.

  \begin{conj}
    \label{conj: crystalline lift implies explicit crystalline lift}
    Let $K/\Qp$ be a finite totally ramified extension, and let
    $\rhobar:G_K\to\GL_2(\Fpbar)$ be a continuous representation. Let
    $a\in\Z^2_+$ be a Serre weight, and suppose
    that for some lift $\lambda\in(\Z^2_+)^{\Hom_{\Qp}(K,\Qpbar)}$, there is
    a continuous crystalline representation
    $\rho:G_{K}\to\GL_2(\Qpbar)$ lifting $\rhobar$, such
    that $\rho$ has Hodge type $\lambda$.

Then $a\in\Wconj(\rbar)$.
  \end{conj}
  We do not know how to prove this conjecture, and we do not directly
  address the conjecture in the rest of this paper. Instead, we
  proceed more indirectly. Proposition \ref{prop: modular of some
    weight implies crystalline lifts exist} is a simple consequence of
  lifting automorphic forms of weight $a$ to forms of weight
  $\lambda$; we may also obtain non-trivial information by lifting to
  forms of weight $0$ and non-trivial type. In this paper, we will
  always consider principal series types. Recall that if $K/\Qp$ is a finite extension the
  \emph{inertial type} of a potentially semistable Galois
  representation $\rho:G_K\to\GL_n(\Qpbar)$ is the restriction to
  $I_K$ of the corresponding Weil-Deligne representation. In this
  paper we normalise this definition as in the appendix to
  \cite{MR1639612}, so that for example the inertial type of a finite
  order character is just the restriction to inertia of that
  character.

\begin{prop}
  \label{prop: modular of some weight implies potentially BT lifts
    exist}Let $F$ be an imaginary CM field with maximal totally real
  subfield $F^+$, and suppose that $F/F^+$ is unramified at all finite
  places, that every place of $F^+$ dividing $p$ has residue field
  $\Fp$ and splits completely in
  $F$, and that $[F^+:\Q]$ is even. Suppose that $p>2$, and that
  $\rbar:G_F\to\GL_2(\Fpbar)$ is an irreducible modular representation
 with split ramification. Let $a\in(\Z^2_+)_0^S$ be a
  Serre weight. If $\rbar$ is modular of weight $a$, then for each
  place $w|p$ of $F$, there is a continuous potentially semistable
  representation $\rho_w:G_{F_w}\to\GL_2(\Qpbar)$ lifting
  $\rbar|_{G_{F_w}}$, such that $\rho_w$ has Hodge type $0$ and
  inertial type $\omegat^{a_1}\oplus\omegat^{a_2}$. (Here $\omegat$ is
  the Teichm\"uller lift of $\omega$.)  Furthermore, $\rho_w$ is
  potentially crystalline unless $a_{1}-a_{2}=p-1$ and  $\rbar|_{G_{F_w}}\cong
  \begin{pmatrix}
    \chibar\epsilonbar&*\\0&\chibar
  \end{pmatrix}
$ for some character $\chibar$.
\end{prop}
\begin{proof}
  This may be proved in exactly the same way as Lemma 3.4 of
  \cite{geesavitttotallyramified}, working in the setting of
  \cite{blggU2} (cf. the proof of Lemma 3.1.1 of \cite{blggU2}). Note
  that if $\rho_w$ is not potentially crystalline, then it is
  necessarily a twist of an extension of the trivial character by the
  cyclotomic character.
\end{proof}

\section{Realising local representations globally}\label{sec:local to
  global}\subsection{}We now recall a result from the forthcoming paper \cite{geekisin}
which allows us to realise local representations globally, in order to
apply the results of Section~\ref{ss:global} in a purely local
setting.

\begin{thm}
  \label{thm: the final local-to-global result} Suppose that $p>2$,
  that $K/\Qp$ is a finite extension, and let
  $\rbar_K:G_K\to\GL_2(\Fpbar)$ be a continuous representation. Then
  there is an imaginary CM field $F$ and a continuous irreducible
  representation $\rbar:G_F\to\GL_2(\Fpbar)$ such that, if $F^+$ denotes the maximal totally real subfield of $F$,
  \begin{itemize}
  \item each place $v|p$ of $F^+$ splits in $F$ and has $F^+_v\cong
    K$,
  \item for each place $v|p$ of $F^+$, there is a place $\tv$ of $F$
    lying over $F^+$ with $\rbar|_{G_{F_\tv}}$ isomorphic to an
    unramified twist of $\rbar_K$,
  \item $\zeta_p\notin F$,
  \item $\rbar$ is unramified outside of $p$,
  \item $\rbar$ is modular in the sense of \cite{blggU2}, and
  \item $\rbar(G_{F(\zeta_p)})$ is adequate.
  \end{itemize}
\end{thm}

\begin{proof}We sketch the proof; the full details will appear in
  \cite{geekisin}. The argument is a straightforward application of
  potential modularity techniques. First, an application of
  Proposition 3.2 of \cite{frankII} supplies a totally real field $L^+$ and a continuous irreducible
  representation $\rbar:G_{L^+}\to\GL_2(\Fpbar)$ such that
  \begin{itemize}
  \item for each place $v|p$ of $L^+$, $L^+_v\cong K$ and
    $\rbar|_{L^+_v}\cong\rbar_K$,
 \item for each place $v|\infty$ of $L^+$, $\det\rbar(c_v)=-1$, where
    $c_v$ is a complex conjugation at $v$, and
  \item there is a non-trivial finite extension $\F/\F_p$ such that
    $\rbar(G_{L^+})=\GL_2(\F)$.
\end{itemize}
By a further base change one can also arrange that  $\rbar|_{G_{L^+_v}}$ is unramified
at each finite place $v\nmid p$ of $L^+$.

By  Lemma 6.1.6 of \cite{blggord} and the proof of
Proposition 7.8.1 of \cite{0905.4266}, $\rbar_K$ admits a potentially
Barsotti-Tate lift, and one may then apply  Proposition 8.2.1 of
\cite{0905.4266} to deduce that there is a finite totally real Galois
  extension $F^+/L^+$ in which all primes of $L^+$ above $p$ split
  completely, such that $\rbar|_{G_{F^+}}$ is modular in the sense
  that it is congruent to the Galois representation associated to some
  Hilbert modular form of parallel weight $2$.

By the theory of base change between $\GL_2$ and unitary groups
(\textit{cf.} section 2 of \cite{blggU2}), it now suffices to show that
there is a totally imaginary quadratic extension $F/F^+$ and a
character $\thetabar:G_F\to\Fpbar^\times$ such that
$\rbar|_{G_F}\otimes\thetabar$ has multiplier~$\epsilonbar^{-1}$ and
such that for each place $v|p$ of $F^+$, there is a place $\tv$ of $F$
lying over $v$ with $\thetabar|_{G_{F_{\tv}}}$ unramified. The
existence of such a character is a straightforward exercise in class
field theory, and follows for example from Lemma 4.1.5 of \cite{cht}.
\end{proof}

\section{Congruences}\label{sec: congruences to weight 0}\subsection{} Having realised a local mod $p$
representation globally, we can now use the results explained in
Section \ref{sec: serre
  weight definitions} to deduce non-trivial local consequences.
\begin{thm}
  \label{thm: explicit weight implies pot BT lift}Let $p>2$ be prime,
  let $K/\Qp$ be a finite totally ramified extension, and let
  $\rhobar:G_K\to\GL_2(\Fpbar)$ be a continuous representation. Let
  $a\in\Wconj(\rhobar)$ be a Serre weight. Then there is a continuous
  potentially semistable representation $\rho:G_K\to\GL_2(\Qpbar)$
  lifting $\rhobar$, such that $\rho$ has Hodge type $0$ and inertial
  type $\omegat^{a_1}\oplus\omegat^{a_2}$. Furthermore, $\rho$ is
  potentially crystalline unless $a_{1}-a_{2}=p-1$ and $\rhobar\cong
  \begin{pmatrix}
    \chibar\epsilonbar&*\\0&\chibar
  \end{pmatrix}$ for some character $\chibar$.

\end{thm}
\begin{proof}  By Theorem \ref{thm: the final local-to-global result}, there is
  an imaginary CM field $F$ and a modular representation
  $\rbar:G_F\to\GL_2(\Fpbar)$ such that
  \begin{itemize}
  \item for each place $v|p$ of $F^+$, $v$ splits in $F$ as
    $\tv\tv^c$, and we have $F_\tv\cong K$, and $\rbar|_{G_{F_\tv}}$ is
    isomorphic to an unramified twist of $\rhobar$,
\item $\rbar$ is unramified outside of $p$,
\item $\zeta_p\notin F$, and
\item $\rbar(G_{F(\zeta_p)})$ is adequate.
\end{itemize}Now, since the truth of the result to be proved is
obviously unaffected by making an unramified twist (if $\rhobar$ is
replaced by a twist by an unramified character $\overline{\theta}$, one may
replace $\rho$ by a twist by an unramified
lift of $\overline{\theta}$), we may without loss of
generality suppose that $\rbar|_{G_{F_w}}\cong\rhobar$. Let
$b\in(\Z^2_+)_0^{S}$ be the Serre weight such that
 $b_\tv=a$ for each place $v|p$ of $F^+$, where $S$ denotes the set of
 places of $F$ above $p$. By Remark \ref{rem: conjectured weights independent of
      unramified twist}, $b\in\Wconj(\rbar)$. Then by Theorem \ref{thm: explicit local lifts implies Serre
    weight}, $\rbar$ is modular of weight $b$. The result now follows
  from Proposition \ref{prop: modular of some weight implies potentially BT lifts
    exist}.
\end{proof}

\subsection{Spaces of crystalline extensions}\label{subsec: H^1_f}We
now specialise to the setting of Definition \ref{defn: W? niveau
  1}. As usual, we let $K/\Qp$ be a finite totally ramified extension with residue
field $k=\Fp$, ramification index $e$, and uniformiser $\pi$. We fix a Serre weight $a\in\Z^2_+$. We fix a
continuous representation $\rhobar:G_K\to\GL_2(\Fpbar)$, and we assume
that there is:
\begin{itemize}
\item  a decomposition
  $\Hom(\Fp,\Fpbar)=J\coprod J^c$, and
\item an integer
  $0\le \delta\le e-1$ such that \[\rhobar|_{I_K}\cong
  \begin{pmatrix}
 \omega^\delta\prod_{\sigma\in
      J}\omega_{\sigma}^{a_{1}+1}\prod_{\sigma\in
      J^c}\omega_\sigma^{a_{2}}&*\\ 0& \omega^{e-1-\delta}\prod_{\sigma\in
      J^c}\omega_\sigma^{a_{1}+1}\prod_{\sigma\in
      J}\omega_\sigma^{a_{2}}.   \end{pmatrix}\]
\end{itemize}
Note that in general there might be several choices
of $J$, $\delta$. Fix such a choice for the
moment. Consider pairs of characters $\chi_1$,
$\chi_2:G_K\to\Qpbar^\times$ with the properties that:
\begin{enumerate}
\item $\rhobar\cong
  \begin{pmatrix}
    \chibar_1&*\\0&\chibar_2
  \end{pmatrix}$,
\item $\chi_1$ and $\chi_2$ are crystalline, and
\item if we let $S$ denote the set of
  $\Hom_{\Qp}(K,\Qpbar)$, then either
  \begin{enumerate}[(i)]
  \item $J$ is non-empty, and there is one embedding $\tau\in
    S$ with $\HT_\tau(\chi_1)=a_{1}+1$ and
    $\HT_\tau(\chi_2)=a_{2}$, there are $\delta$ embeddings
    $\tau\in S$ with  $\HT_\tau(\chi_1)=1$ and
    $\HT_\tau(\chi_2)=0$, and for the remaining $e-1-\delta$
    embeddings $\tau\in S$ we have  $\HT_\tau(\chi_1)=0$ and
    $\HT_\tau(\chi_2)=1$, or
   \item $J=\emptyset$, and there is one embedding $\tau\in
    S$ with $\HT_\tau(\chi_1)=a_{2}$ and
    $\HT_\tau(\chi_2)=a_{1}+1$, there are $\delta$ embeddings
    $\tau\in S$ with  $\HT_\tau(\chi_1)=1$ and
    $\HT_\tau(\chi_2)=0$, and for the remaining $e-1-\delta$
    embeddings $\tau\in S$ we have  $\HT_\tau(\chi_1)=0$ and
    $\HT_\tau(\chi_2)=1$.
  \end{enumerate}
\end{enumerate}
Note that these properties do not specify the characters $\chi_1$ and
$\chi_2$ uniquely, even in the unramified case, as one is always free
to twist either character by an unramified character which is trivial
mod $p$.  We point out that the Hodge type of any de Rham extension of
$\chi_2$ by $\chi_1$
will be a lift of $a$.  Conversely, by Lemma~6.2 of \cite{geesavitttotallyramified} any $\chi_1,\chi_2$ satisfying
(1) and (2) such that the Hodge type of $\chi_1 \oplus \chi_2$ is a
lift of $a$ will satisfy (3) for a valid choice of $J$ and $\delta$
(unique unless $a=0$).

Suppose now that we have fixed two such characters $\chi_1$ and
$\chi_2$, and we now allow the (line corresponding to the) extension
class of $\rhobar$ in $\Ext_{G_K}(\chibar_2,\chibar_1)$ to vary. We
naturally identify $\Ext_{G_K}(\chibar_2,\chibar_1)$ with
$H^1(G_K,\chibar_1 \chibar_2^{-1})$ from now on.
\begin{defn}
  Let $L_{\chi_1,\chi_2}$ be the subset of
  $H^1(G_K,\chibar_1\chibar_2^{-1})$ such that the corresponding
  representation $\rhobar$ has a crystalline lift $\rho$ of the
  form \[
  \begin{pmatrix}
    \chi_1&*\\0&\chi_2
  \end{pmatrix}.\]
\end{defn}
We have the following variant of Lemma 3.12 of \cite{bdj}.
\begin{lem}
  \label{lem: dimension of H^1_f spaces} $L_{\chi_1,\chi_2}$ is an
  $\Fpbar$-vector subspace of $ H^1(G_K,\chibar_1\chibar_2^{-1})$ of
  dimension $|J|+\delta$, unless
  $\chibar_1=\chibar_2$, in which case it has dimension
  $|J|+\delta+1$.
\end{lem}
\begin{proof} Let $\chi=\chi_1\chi_2^{-1}$.
  Recall that
  $H^1_f(G_K,\Zpbar(\chi))$ is the preimage of
  $H^1_f(G_K,\Qpbar(\chi))$ under the natural map
  $\eta : H^1(G_K,\Zpbar(\chi))\to H^1(G_K,\Qpbar(\chi))$, so that
  $L_{\chi_1,\chi_2}$ is the image of $H^1_f(G_K,\Zpbar(\chi))$ in
  $H^1(G_K,\chibar)$.  The kernel of $\eta$ is precisely the torsion
  part of $H^1_f(G_K,\Zpbar(\chi))$, which (since $\chi\neq 1$,
  e.g. by examining Hodge-Tate weights) is non-zero  if and only if
  $\chibar=1$, in which case it has the form $\kappa^{-1} \Zpbar/\Zpbar$
  for some  $\kappa \in \mf{m}_{\Zpbar}$.

By
Proposition 1.24(2) of \cite{nekovar} we see that $\dim_\Qpbar
H^1_f(G_K,\Qpbar(\chi))=|J|+\delta$, again using $\chi \neq 1$. Since
$H^1(G_K,\Zpbar(\chi))$ is a finitely generated $\Zpbar$-module,
the result follows.
  \end{proof}
  \begin{defn}
    \label{defn: union of H^1_f subspaces}If $\chibar_1$ and
      $\chibar_2$ are fixed, we define $\Lcrys$ to be the subset of
      $H^1(G_K,\chibar_1 \chibar_2^{-1})$ given by the union of the $L_{\chi_1,\chi_2}$
      over all $\chi_1$ and $\chi_2$ as above.
  \end{defn}
Note that $\Lcrys$ is a union of subspaces of possibly varying
dimensions, and as such it is not clear that $\Lcrys$ is itself a
subspace. Note also that the representations $\rhobar$ corresponding
to elements of $\Lcrys$ are by definition precisely those for which
$F_a\in\Wconj(\rhobar)$.

\begin{defn}
  \label{defn: H^1_flat subspace}Let $\Lflat$ be the subset
  of $H^1(G_K,\chibar_1\chibar_2^{-1})$ consisting of classes with the property that
  if $\rhobar\cong
  \begin{pmatrix}
    \chibar_1&*\\0&\chibar_2
  \end{pmatrix}$ is the corresponding representation, then there is a
  finite field $k_E \subset \Fpbar$  and a
  finite flat $k_E$-vector space scheme over $\cO_{K(\pi^{1/(p-1)})}$ with
  generic fibre
  descent data to $K$ of the
  form $ \omega^{a_{1}}\oplus\omega^{a_{2}}$
(see Definition~\ref{defn:dd-of-the-form}) whose generic fibre is $\rhobar$.
\end{defn}
\begin{thm}
  \label{thm: crystalline extension implies flat}Provided that
  $a_{1}-a_{2}\ne p-1$ or that $\chibar_1\chibar_2^{-1}\ne \epsilonbar$,
  $\Lcrys\subset\Lflat$.
\end{thm}
\begin{proof}
   Take a class in $\Lcrys$, and consider the corresponding
  representation  $\rhobar\cong
  \begin{pmatrix}
    \chibar_1&*\\0&\chibar_2
  \end{pmatrix}$. As remarked above, $F_a\in\Wconj(\rhobar)$, so by
  Theorem \ref{thm: explicit weight implies pot BT lift}, $\rhobar$
  has a crystalline lift of Hodge type $0$ and inertial
  type \[\omegat^{a_{1}}\oplus\omegat^{a_{2}},\] and this
  representation can be taken to have coefficients in the ring of
  integers $\cO_E$ of a finite
  extension $E/\Qp$.  Let $\varpi$ be a uniformiser of $\cO_E$, and $k_E$
  the residue field.  Such a representation
  corresponds to a $p$-divisible $\cO_E$-module with generic fibre descent data, and
  taking the $\varpi$-torsion
  gives a finite flat $k_E$-vector space scheme with generic fibre descent
  data whose generic fibre is $\rhobar$.  By  Corollary 5.2 of \cite{geesavittquaternionalgebras} this
descent data has the form $\omega^{a_1} \oplus \omega^{a_2}$.
\end{proof}
In the next section we will make calculations with finite flat group
schemes in order to relate $\Lflat$ and $\Lcrys$.

\section{Finite flat models}\label{sec: finite flat
  models}\subsection{}We work throughout this section in the following setting:
\begin{itemize}
\item $K/\Qp$ is a finite extension with ramification index $e$,
  inertial degree $1$, ring
  of integers $\cO_K$, uniformiser $\pi$ and residue field $\Fp$.
\item $\chibar_1$, $\chibar_2$ are characters
  $G_K\to\Fpbar^\times$.
\item $a\in\Z^2_+$ is a Serre weight.
\item There is a decomposition $\Hom(\Fp,\Fpbar)=J\coprod J^c$, and an integer $0\le
  \delta\le e-1$ such that \[\chibar_1|_{I_K}=\omega^\delta\prod_{\sigma\in
    J}\omega^{a_{1}+1}\prod_{\sigma\in
    J^c}\omega^{a_{2}},\] \[\chibar_2|_{I_K}=\omega^{e-1-\delta}\prod_{\sigma\in
    J^c}\omega^{a_{1}+1}\prod_{\sigma\in
    J}\omega^{a_{2}}.\]
\end{itemize}
Note in particular that $(\chibar_1\chibar_2)|_{I_K}=\omega^{a_1+a_2+e}$.

Let $K_1:=K(\pi^{1/(p-1)})$. Let $k_E$ be a finite extension of $\Fp$
such that $\chibar_1,\chibar_2$ are defined over $k_E$; for the moment
$k_E$ will be fixed, but eventually it will be allowed to vary.
We wish to consider the representations  $\rhobar\cong
  \begin{pmatrix}
    \chibar_1&*\\0&\chibar_2
  \end{pmatrix}$ such that there is a finite flat $k_E$-vector space
  scheme $\cG$ over $\cO_{K_1}$ with generic fibre descent data to $K$ of the form
  $\omega^{a_1}\oplus\omega^{a_2}$ (see Definition~\ref{defn:dd-of-the-form}), whose generic fibre is
  $\rhobar$.

In order to do, we will work with Breuil modules with descent data
from $K_1$ to~$K$. We
recall the necessary definitions from
\cite{geesavittquaternionalgebras}.

Fix $\pi_1$, a $(p-1)$-st root of $\pi$ in $K_1$. Write
$e'=e(p-1)$. The category $\BrMod_{\dd}$
consists of quadruples $(\mathcal{M},\Fil^1
\mathcal{M},\phi_{1},\{\widehat{g}\})$ where:

\begin{itemize}\item $\mathcal{M}$ is a finitely generated free
  $k_E[u]/u^{e'p}$-module,
\item $\Fil^1 \M$ is a $k_E[u]/u^{e'p}$-submodule of $\M$ containing $u^{e'}\M$,
\item $\phi_{1}:\Fil^1\M\to\M$ is $k_E$-linear and $\phi$-semilinear
  (where $\phi:\Fp[u]/u^{e'p}\to \Fp[u]/u^{e'p}$ is the $p$-th power map)
  with image generating $\M$ as a $k_E[u]/u^{e'p}$-module, and
\item $\widehat{g}:\M\to\M$ for each $g\in\Gal(K_1/K)$ are additive
  bijections that preserve $\Fil^1 \M$, commute with the $\phi_1$-,
  and $k_E$-actions, and satisfy $\widehat{g}_1\circ
  \widehat{g}_2=\widehat{g_1\circ g}_2$ for all
  $g_1,g_2\in\Gal(K_1/K)$; furthermore $\widehat{1}$ is the identity,
  and if $a\in k_E$, $m\in\M$ then
  $\widehat{g}(au^{i}m)=a((g(\pi)/\pi)^{i})u^{i}\widehat{g}(m)$.\end{itemize}

The category $\BrMod_{\dd}$ is equivalent to the category of finite
flat $k_E$-vector space schemes over $\mathcal{O}_{K_1}$ together with
descent data on the generic fibre from $K_1$ to~$K$
(this equivalence depends on $\pi_1$); see \cite{sav06}, for instance. We obtain the associated
$G_{K}$-representation (which we will refer to as the generic fibre)
of an object of $\BrMod_{\dd,K_1}$ via the covariant functor
$T_{\st,2}^{K}$ (which is defined immediately before Lemma 4.9 of
\cite{MR2137952}).

\begin{defn}
  \label{defn:dd-of-the-form}
  Let $\M$ be an object of $\BrMod_{\dd}$ such that the underlying
  $k_E$-module has rank two. We say that the finite flat $k_E$-vector
  space scheme corresponding to $\M$ \emph{has descent data
  of the form} $\omega^{a_1} \oplus \omega^{a_2}$ if $\M$ has a basis
 $e_1,e_2$ such that $\widehat{g}(e_i) = \omega^{a_i}(g) e_i$. (Here
 we abuse notation by identifying an element of $G_K$ with its image
 in $\Gal(K_1/K)$.)
\end{defn}

We now consider a finite flat group scheme with generic fibre descent data $\cG$ as above. By a standard scheme-theoretic
closure argument, $\chibar_1$ corresponds to a finite flat subgroup
scheme with generic fibre descent data
$\mathcal{H}$ of $\mathcal{G}$, so we begin by analysing the possible
finite flat group schemes corresponding to characters.

Suppose now that $\M$ is an object of $\BrMod_{\dd}$. The rank
one objects of $\BrMod_{\dd}$ are classified as follows.

\begin{prop} \label{prop:rank one breuil modules} With our fixed choice of uniformiser
 $\pi$, every rank one object of $\BrMod_{\dd}$ has the form:
\begin{itemize}
\item $\M = (k_E[u]/u^{e'p}) \cdot v $,
\item $\Fil^1 \M = u^{x(p-1)} \M$,
\item $\phi_1( u^{x(p-1)} v) = cv$ for some $c \in k_E^{\times}$, and
\item $\widehat{g}(v) = \omega(g)^kv$ for all $g \in \Gal(K_1/K)$,
\end{itemize}
where $0 \le x \le e$ and $0 \le k< p-1$ are
integers.

Then $T_{\st,2}^{K}(\M) =
\omega^{k + x} \cdot \ur_{c^{-1}}$, where $\ur_{c^{-1}}$ is the
unramified character taking an arithmetic Frobenius element to
$c^{-1}$.\end{prop}
\begin{proof}
  This is a special case of Proposition 4.2 and Corollary 4.3 of
  \cite{geesavittquaternionalgebras}.
\end{proof}

Let $\M$ (or $\M(x)$) be the rank one Breuil
module with $k_E$-coefficients and
descent data from $K_1$ to $K$ corresponding to $\mathcal{H}$, and
write $\M$ in the form given by Proposition \ref{prop:rank one breuil
  modules}. Since $\mathcal{G}$ has descent data of the form
$\omega^{a_1}\oplus\omega^{a_2}$,
we must have  $\omega^k \in \{\omega^{a_1},\omega^{a_2}\}$.

\subsection{Extensions} Having determined the rank one characters, we
now go further and compute the possible extension
classes. By a scheme-theoretic closure argument, the Breuil module
$\cP$ corresponding to $\cG$ is an extension of $\cN$ by
$\cM$, where $\cM$ is as in the previous section, and $\cN$ (or
$\cN(y)$)  is defined
by \begin{itemize}
\item $\N = (k_E[u]/u^{e'p}) \cdot w $,
\item $\Fil^1 \N = u^{y(p-1)} \N$,
\item $\phi_1( u^{y(p-1)} v) = dw$ for some $d \in k_E^{\times}$, and
\item $\widehat{g}(v) = \omega(g)^lv$ for all $g \in \Gal(K_1/K)$,
\end{itemize}
where $0 \le y \le e$ and $0 \le l< p-1$ are
integers. Now, as noted above, the descent data for $\mathcal{G}$ is of the form
$\omega^{a_1}\oplus\omega^{a_2}$, so we must have that either $\omega^k=\omega^{a_1}$
and $\omega^l=\omega^{a_2}$, or $\omega^{k}=\omega^{a_2}$ and $\omega^l=\omega^{a_1}$. Since by definition we have
$(\chibar_1\chibar_2)|_{I_K}=\omega^{a_1+a_2+e}$, we see from
Proposition \ref{prop:rank one breuil modules} that \[x+y\equiv e\pmod{p-1}.\]

  We have the following classification of extensions of $\cN$ by $\cM$.
\begin{prop}\label{prop: possible extensions of Breuil modules} Every extension of $\cN$ by
  $\cM$ is isomorphic to exactly one of the form
  \begin{itemize}
\item $\cP = (k_E[u]/u^{e'p}) \cdot v + (k_E[u]/u^{e'p}) \cdot w $,
\item $\Fil^1 \cP =(k_E[u]/u^{e'p}) \cdot u^{x(p-1)} v +
  (k_E[u]/u^{e'p}) \cdot (u^{y(p-1)}w+\lambda v) $,
\item $\phi_1(u^{x(p-1)} v) = cv$, $\phi_1(u^{y(p-1)}w+\lambda v)=dw$,
\item $\widehat{g}(v) =\omega^k(g)v$ and  $\widehat{g}(w) =\omega^l(g)w$ for all $g \in \Gal(K_1/K)$,
\end{itemize}where $\lambda\in u^{\max\{0,(x+y-e)(p-1)\}}k_E[u]/u^{e'p}$
has all nonzero terms of degree congruent to $l-k$ modulo $p-1$, and has all terms
of degree less than $x(p-1)$, unless $\chibar_1=\chibar_2$ and $x\ge y$,
in which case it may additionally have a term of degree $px-y$.
\end{prop}
\begin{proof}
  This is a special case of Theorem 7.5 of \cite{MR2004122}, with the
  addition of $k_E$-coefficients in place of $\Fp$-coefficients. When
  $K$ (in the notation of \emph{loc.~cit.}) is totally ramified over $\Qp$, the
  proof of \emph{loc.~cit.} is argued in precisely the same manner when
  coefficients are added, taking care to note the following changes:
\begin{itemize}
\item Replace Lemma 7.1 of \emph{loc.~cit.} (i.e., Lemma 5.2.2 of
  \cite{MR1839918}) with Lemma 5.2.4 of \cite{MR1839918}  (with
  $k'=k_E$ and $k=\Fp$ in the notation of that Lemma).  In particular
  replace $t^l$ with $\phi(t)$ wherever it appears in the proof, where~$\phi$
is the $k_E$-linear endomorphism of $k_E[u]/u^{e'p}$ sending
  $u^i$ to $u^{pi}$.

\item Instead of applying Lemma 4.1 of \cite{MR2004122}, note that the
 cohomology group
 $H^1(\Gal(K_1/K),k_E[u]/u^{e'p})$ vanishes because $\Gal(K_1/K)$ has prime-to-$p$
 order while $k_E[u]/u^{e'p}$ has $p$-power order.

\item Every occurrence of $T_{i}^l$ in the proof (for any subscript $i$) should be replaced with
  $T_{i}$. In the notation of \cite{MR2004122} the element $\eta$ is
  defined when the map $\alpha \mapsto  (1-b/a)\alpha$ on $k_E$ is not
  surjective, i.e., when $a=b$; we may then take $\eta=1$.

\item The coefficients of $h,t$ are permitted to lie in $k_E$
  (i.e., they are not constrained to lie in any particular proper subfield).
\end{itemize}
\end{proof}

Note that the recipe for $\cP$ in the statement of
Proposition~\ref{prop: possible extensions of Breuil modules} defines
an extension of $\cN$ by $\cM$ provided that $\lambda$ lies in  $u^{\max\{0,(x+y-e)(p-1)\}}k_E[u]/u^{e'p}$
and has all nonzero terms of degree congruent to $l-k$ modulo $p-1$
(\emph{cf.} the discussion in Section 7 of \cite{MR2004122}).  Denote
this Breuil module by $\cP(x,y,\lambda)$.    Note that $c$ is fixed
while $x$ determines
$k$, since we require $\omega^{k+x} \cdot \ur_{c^{-1}} =
\chibar_1$; similarly $d$ is fixed and $y$ determines $l$.  So this notation
is reasonable.

We would like to compare the generic fibres of extensions of different
choices of $\cM$ and $\cN$. To this end, we have the following
result.  Write
$\chibar_1|_{I_K}=\omega^\alpha$, $\chibar_2|_{I_K}=\omega^\beta$.
\begin{prop}
  \label{prop: comparing extensions}The Breuil module $\cP(x,y,\lambda)$ has the same generic fibre as the Breuil module $\cP'$,
  where  \begin{itemize}
\item $\cP' = (k_E[u]/u^{e'p}) \cdot v' + (k_E[u]/u^{e'p}) \cdot w' $,
\item $\Fil^1 \cP' =(k_E[u]/u^{e'p}) \cdot u^{e(p-1)} v' +
  (k_E[u]/u^{e'p}) \cdot (w'+u^{p(e-x)+y}\lambda v') $,
\item $\phi_1(u^{e(p-1)} v') = cv'$, $\phi_1(w'+u^{p(e-x)+y}\lambda v')=dw'$,
\item $\widehat{g}(v') =\omega^{\alpha-e}(g)v'$ and  $\widehat{g}(w') =\omega^{\beta}(g)w'$ for all $g \in \Gal(K_1/K)$.
\end{itemize}
\end{prop}
\begin{proof}
  Consider the Breuil module $\cP''$ defined by  \begin{itemize}
\item $\cP'' = (k_E[u]/u^{e'p}) \cdot v'' + (k_E[u]/u^{e'p}) \cdot w'' $,
\item $\Fil^1 \cP'' =(k_E[u]/u^{e'p}) \cdot u^{e(p-1)} v'' +
  (k_E[u]/u^{e'p}) \cdot (u^{y(p-1)}w''+u^{p(e-x)}\lambda v'') $,
\item $\phi_1(u^{e(p-1)} v'') = cv''$, $\phi_1(u^{y(p-1)}w''+u^{p(e-x)+y}\lambda v'')=dw''$,
\item $\widehat{g}(v'') =\omega^{k+x-e}(g)v''$ and  $\widehat{g}(w'') =\omega^{l}(g)w''$ for all $g \in \Gal(K_1/K)$.
\end{itemize}
(One checks without difficulty that this \emph{is} a Breuil module.  For instance the condition
on the minimum degree of terms appearing in $\lambda$ guarantees that
$\Fil^1 \cP''$ contains $u ^{e'}\cP''$.)  Note that $k+x\equiv \alpha\pmod{p-1}$,
  $l+y\equiv\beta\pmod{p-1}$.   We claim that $\cP$, $\cP'$ and $\cP''$ all have the
same generic fibre. To see this, one can check directly that there is a morphism
$\cP\to\cP''$ given by \[v\mapsto u^{p(e-x)}v'',\ w\mapsto w'',\]and a
morphism $\cP'\to\cP''$ given by \[v'\mapsto v'',\ w'\mapsto
u^{py}w''.\] By Proposition 8.3 of \cite{MR2004122}, it is enough to
check that the kernels of these maps do not contain any free
$k_E[u]/(u^{e'p})$-submodules, which is an immediate consequence of
the inequalities $p(e-x),py<e'p$.
\end{proof}

\begin{rem}
  \label{rem:extension-classes}
  We note for future reference that while the classes in
  $H^1(G_K,\chibar_1 \chibar_2^{-1})$ realised by $\cP(x,y,\lambda)$ and
  $\cP'$ may not coincide, they differ at most by multiplication
  by a $k_E$-scalar.  To see this, observe that the maps $\cP \to
  \cP''$ and $\cP' \to \cP''$ induce $k_E$-isomorphisms on the
  rank one sub- and quotient Breuil modules.
\end{rem}

We review the constraints on the integers $x,y$: they must lie
between $0$ and~$e$, and if we let $k,l$ be the residues of
$\alpha-x,\beta-y \pmod{p-1}$ in the interval $[0,p-1)$ then we must
have $\{\omega^k,\omega^l\} = \{\omega^{a_1},\omega^{a_2}\}$. Call such a pair $x,y$ \emph{valid}.
Note that $l-k \equiv \beta-\alpha + x - y \pmod{p-1}$ for any valid pair.

\begin{cor}
  \label{cor:comparison-of-good-models}
Let $x',y'$ be another valid pair.
  Suppose that $x' + y' \le e$ and $p(x'-x)+(y -y') \ge 0$.   Then $\cP(x,y,\lambda)$ has
  the same generic fibre as
  $\cP(x',y',\lambda')$, where $\lambda' = u^{p(x'-x)+(y-y')} \lambda$.
\end{cor}

\begin{proof}
 The Breuil module $\cP(x',y',\lambda')$ is well-defined: one checks
 from the definition that the congruence
 condition on the degrees of the nonzero terms in $\lambda'$ is
 satisfied, while since $x'+y' \le e$
 there is no condition on the lowest degrees appearing in $\lambda'$.
 Now the result is immediate from Proposition~\ref{prop: comparing extensions},
  since $u^{p(e-x)+y}\lambda = u^{p(e-x')+y'}\lambda'$.
\end{proof}

Recall that $x+y \equiv e\pmod{p-1}$, so that $x$ and $e-y$ have the
same residue modulo $p-1$.   It follows that if $x,y$ is a valid pair
of parameters, then so is $e-y,y$.   Let $X$ be the largest
value of $x$ over all valid pairs $x,y$, and similarly $Y$ the smallest value of $y$;
then $Y=e-X$, since if if we had $Y > e-X$ then $e-Y$ would be a
smaller  possible value for $x$.

\begin{cor}
  \label{cor:generic-fibres-all-occur-extremally}
  The module $\cP(x,y,\lambda)$ has the same generic fibre as
  $\cP(X,Y,\mu)$ where $\mu \in k_E[u]/u^{e'p}$ has all nonzero terms of degree congruent to $\beta-\alpha+X-Y$ modulo $p-1$, and has all terms
of degree less than $X(p-1)$, unless $\chibar_1=\chibar_2$,
in which case it may additionally have a term of degree
$pX-Y$.
\end{cor}

\begin{proof}
 Since $X+Y=e$ and $p(X-x)+(y-Y) \ge 0$ from the choice of $X,Y$, the
 previous Corollary shows that $\cP(x,y,\lambda)$ has the same generic
 fibre as some $\cP(X,Y,\lambda')$; by Proposition~\ref{prop:
   possible extensions of Breuil modules}   this has the same generic
 fibre as $\cP(X,Y,\mu)$  for $\mu$ as in the statement.  (Note that
 if $\chibar_1=\chibar_2$ then automatically $X \ge Y$, because in this
 case if $x,y$ is a valid pair then so is $y,x$.)
\end{proof}

\begin{prop}
  \label{prop:computation of the dimension of Lflat}Let $X$ be as
  above, i.e., $X$ is the maximal integer
  such
  that
  \begin{itemize}
  \item $0\le X\le e$, and
  \item either $\chibar_1|_{I_K}=\omega^{a_1+X}$ or
    $\chibar_1|_{I_K}=\omega^{a_2+X}$.
  \end{itemize}
Then $\Lflat$is an $\Fpbar$-vector space of dimension at most
$X$, unless $\chibar_1=\chibar_2$, in which case it has
dimension at most $X+1$.
\end{prop}
\begin{proof}
  Let $L_{\mathrm{flat},k_E} \subset \Lflat$ consist of the  classes $\eta$ such that the containment $\eta \in \Lflat$ is
  witnessed by a $k_E$-vector space scheme with generic fibre descent
  data.  By
  Corollary~\ref{cor:generic-fibres-all-occur-extremally} and Remark~\ref{rem:extension-classes} these are exactly
  the classes arising from the Breuil modules $\cP(X,Y,\mu)$ with
  $k_E$-coefficients as in
  Corollary~\ref{cor:generic-fibres-all-occur-extremally}.  These
  classes form a $k_E$-vector space (since they are \emph{all} the
  extension classes arising from extensions of $\cN(Y)$ by $\cM(X)$),
  and by counting the (finite) number of possibilities for $\mu$ we see
  that $\dim_{k_E} L_{\mathrm{flat},k_E}$ is at most $X$ (resp.
  $X+1$ when $\chibar_1=\chibar_2$).

  Since $L_{\mathrm{flat},k_E}  \subset L_{\mathrm{flat},k'_E}$ if
  $k_E \subset k'_E$ it follows easily that $\Lflat = \cup_{k_E}
  L_{\mathrm{flat},k_E}$ is an $\Fpbar$-vector space of dimension at
  most $X$ (resp. $X+1$).
\end{proof}
We can now prove our main local result, the promised relation between $\Lflat$ and
$\Lcrys$. \begin{thm}
  \label{thm: crystalline equals flat}Provided that either $a_1-a_2\ne
  p-1$ or $\chibar_1\chibar_2^{-1}\ne\epsilonbar$, we
  have $\Lflat=\Lcrys$.
\end{thm}
\begin{proof}By Theorem \ref{thm: crystalline extension implies flat},
  we know that $\Lcrys\subset\Lflat$, so by Proposition~\ref{prop:computation of the dimension of Lflat} it suffices to show that
  $\Lcrys$ contains an $\Fpbar$-subspace of dimension
  $X$ (respectively $X+1$ if $\chibar_1 = \chibar_2$). Since $\Lcrys$ is the union of the spaces
  $L_{\chi_1,\chi_2}$, it suffices to show that one of these spaces
  has the required dimension. Let $X$ be as in the statement of
  Proposition \ref{prop:computation of the dimension of Lflat}, so
  that $X$ is maximal in $[0,e]$ with the property that either  $\chibar_1|_{I_K}=\omega^{a_1+X}$ or
    $\chibar_1|_{I_K}=\omega^{a_2+X}$. Note that by the assumption
    that there is a decomposition
  $\Hom(\Fp,\Fpbar)=J\coprod J^c$, and an integer
  $0\le \delta\le e-1$ such that \[\rhobar|_{I_K}\cong
  \begin{pmatrix}
\omega^\delta \prod_{\sigma\in
      J}\omega_{\sigma}^{a_{1}+1}\prod_{\sigma\in
      J^c}\omega_\sigma^{a_{2}}&*\\ 0& \omega^{e-1-\delta}\prod_{\sigma\in
      J^c}\omega_\sigma^{a_{1}+1}\prod_{\sigma\in
      J}\omega_\sigma^{a_{2}}  \end{pmatrix},\]we see that
  if $X=0$ then  $\chibar_1|_{I_K}=\omega^{a_2}$ (and $J$ must be empty).

  If $\chibar_1|_{I_K}=\omega^{a_2+X}$ then we take $J$ to be empty
  and we take $\delta=X$; otherwise $X > 0$ and $\chibar_1|_{I_K} =
  \omega^{a_1+X}$, and we can take $J^c$ to be empty and
  $\delta=X-1$. In either case, we may define characters $\chi_1$ and
  $\chi_2$ as in Section \ref{subsec: H^1_f}, and we see from Lemma
  \ref{lem: dimension of H^1_f spaces} that
  $\dim_{\Fpbar}L_{\chi_1,\chi_2}=X$ unless $\chibar_1=\chibar_2$, in
  which case it is $X+1$. The result follows.\end{proof}
As a consequence of this result, we can also address the question of
the relationship between the different spaces $L_{\chi_1,\chi_2}$ for
a fixed Serre weight $a\in\Wconj(\rhobar)$. If $e$ is large, then
these spaces do not necessarily have the same dimension, so they
cannot always be equal. However, it is usually the case that the
spaces of maximal dimension coincide, as we can now see.
\begin{cor}
  \label{cor: independence of lift for H^1_f}If either $a_1-a_2\ne
  p-1$ or $\chibar_1\chibar_2^{-1}\ne\epsilonbar$, then
  the spaces $L_{\chi_1,\chi_2}$ of maximal dimension are all equal.
\end{cor}
\begin{proof}
  In this case $\dim_{\Fpbar} L_{\chi_1,\chi_2}=\dim_{\Fpbar}\Lcrys$
  by the proof of Theorem \ref{thm: crystalline equals flat}, so we
  must have $L_{\chi_1,\chi_2}=\Lcrys$.
\end{proof}
Finally, we determine $\Lcrys$ in the one remaining case, where the
spaces $L_{\chi_1,\chi_2}$ of maximal dimension no longer coincide.
\begin{prop}
  \label{prop: Lcrys in the exceptional case}Suppose that
  $a_1-a_2=p-1$ and that $\chibar_1\chibar_2^{-1}=\epsilonbar$. Then $\Lcrys=H^1(G_K,\epsilonbar)$.
\end{prop}
\begin{proof}We prove this in a similar fashion to the proof of Lemma
  6.1.6 of \cite{blggord}.  By twisting we can reduce to the case
  $(a_1,a_2)=(p-1,0)$.  Let $L$ be a given line in
  $H^1(G_K,\epsilonbar)$, and choose an unramified character $\psi$
  with trivial reduction.  Let
  $\chi$ be some fixed crystalline character of $G_K$ with Hodge-Tate weights
  $p,1,\dots,1$ such that $\chibar=\epsilonbar$. Let $E/\Qp$ be a finite extension with ring
  of integers $\cO$, uniformiser $\varpi$ and residue field $\F$, such
  that $\psi$ and $\chi$ are defined over $E$ and $L$ is defined over $\F$. Since any extension of $1$ by $\chi\psi$ is
  automatically crystalline, it suffices to show that we can choose
  $\psi$ so that $L$ lifts to $H^1(G_K,\cO(\psi\chi))$.

Let $H$ be the
  hyperplane in $H^1(G_K,\bb{F})$ which annihilates $L$ under the Tate
  pairing. Let $\delta_1 : H^1(G_K,\bb F(\overline{\epsilon})) \to
  H^2(G_K,\mc{O}(\psi\chi))$ be the map coming from
  the exact sequence $0\to \mc{O}(\psi\chi)\stackrel{\varpi}{\to}\mc
  O(\psi\chi)\to \bb F(\overline{\epsilon})\to 0$ of
  $G_K$-modules. We need to show that $\delta_1(L)=0$ for some choice
  of $\psi$.

  Let $\delta_0$ be the map
  $H^0(G_K,(E/\mc{O})(\psi^{-1}\chi^{-1}\epsilon)) \to
  H^{1}(G_K,\bb{F})$ coming from the exact sequence $0 \to \bb{F} \to
  (E/\mc{O})(\psi^{-1}\chi^{-1}\epsilon) \stackrel{\varpi}{\to}
  (E/\mc{O})(\psi^{-1}\chi^{-1}\epsilon) \to 0$ of $G_K$-modules.  By
  Tate local duality, the condition that $L$ vanishes under the map
  $\delta_1$ is equivalent to the condition that the image of the map
  $\delta_0$ is contained in $H$.  Let $n \geq 1$ be the largest
  integer with the property that $\psi^{-1}\chi^{-1}\epsilon \equiv 1
  \pmod{\varpi^n}$. Then we can write $\psi^{-1}\chi^{-1}\epsilon(x)=
  1+\varpi^n \alpha(x)$ for some function $\alpha : G_K \to
  \mc{O}$. Let $\overline{\alpha}$ denote $\alpha \pmod{\varpi} : G_K
  \to \bb{F}$. Then $\overline{\alpha}$ is additive and the choice of
  $n$ ensures that it is non-trivial. It is straightforward to check
  that the image of the map $\delta_0$ is the line spanned by
  $\overline{\alpha}$. If $\overline{\alpha}$ is in $H$, we are
  done. Suppose this is not the case. We break the rest of the proof
  into two cases.

  \medskip{\sl Case 1: $L$ is
    tr\`es ramifi\'e:}  To begin, we observe that it is
  possible to have chosen
  $\psi$ so that
  $\overline{\alpha}$ is ramified. To see this, let $m$ be the largest integer with the property that
$(\psi^{-1} \chi^{-1} \epsilon)|_{I_K} \equiv 1 \pmod{\varpi^m}$.   Note that $m$ exists since the
Hodge-Tate weights of $\psi^{-1}\chi^{-1}\epsilon$ are not all $0$.
If $m = n$ then we are done, so assume instead that $m >n$. Let $g\in
G_K$ be a lift of $\Frob_K$. We claim that
$\psi^{-1}\chi^{-1}\epsilon(g)= 1 +\varpi^{n} \alpha(g)$ such that
$\alpha (g) \not \equiv 0 \pmod{\varpi}$. In fact, if $\alpha
(g)\equiv 0 \pmod{\varpi}$ then $\psi^{-1}\chi^{-1}\epsilon(g) \in  1
+ \varpi^{n+1} \mc{O}_K$. Since $m > n$ we see that
$\psi^{-1}\chi^{-1}\epsilon(G_K) \subset 1 + \varpi^{n+1} \mc{O}_K$
and this contradicts the selection of $n$. Now define a unramifed
character $\psi'$ with trivial reduction by setting $\psi' (g) =
1 - \varpi^n \alpha (g)$. After replacing $\psi$ by $\psi \psi'$ we
see that $n$ has increased but $m$ has not changed.  After finitely
many iterations of this procedure we have  $m=n$, completing the
claim.

Suppose, then, that $\overline{\alpha}$ is ramified.    The fact that $L$ is tr\`es
  ramifi\'e implies that $H$ does not contain the unramified line in
  $H^1(G_K,\bb{F})$. Thus there is a unique $\overline{x} \in
  \bb{F}^\times$ such that $\overline{\alpha}+u_{\overline{x}} \in H$
  where $u_{\overline{x}}: G_K\to \bb{F}$ is the unramified
  homomorphism sending $\Frob_K$ to $\overline{x}$. Replacing $\psi$ with $\psi$ times
  the unramified character sending $\Frob_K$ to $(1+\varpi^n x)^{-1}$,
  for $x$ a lift of $\overline{x}$, we are done.

  \medskip{\sl Case 2: $L$ is peu ramifi\'e:} Making a ramified
  extension of $\mc{O}$ if necessary, we can and do assume that $n\geq
  2$. The fact that $L$ is peu ramifi\'e implies that $H$ contains the
  unramified line. It follows that if we replace $\psi$ with $\psi$
  times the unramified character sending $\Frob_K$ to $1+\varpi$, then
  we are done (as the new $\overline{\alpha}$ will be unramified).
\end{proof}

\section{Global consequences}\label{sec: global
  consequences}\subsection{}We now deduce our main global results,
using the main theorems of \cite{blggU2} together with our local
results to precisely determine the set of Serre weights for a global
representation in  the totally ramified case.

\begin{prop}
  \label{prop: semisimple elimination if totally ramified}Let $F$ be an imaginary CM field with maximal totally real
  subfield $F^+$, and suppose that $F/F^+$ is unramified at all finite
  places, that every place of $F^+$ dividing $p$ splits completely in
  $F$, and that $[F^+:\Q]$ is even. Suppose that $p>2$, and that
  $\rbar:G_F\to\GL_2(\Fpbar)$ is an irreducible modular representation
  with split ramification. Let
  $a\in(\Z^2_+)_0^S$ be a Serre
  weight such that $\rbar$ is modular of weight $a$. Let $w$ be a
  place of $F$ such that $F_w/\Qp$ is totally ramified of degree $e$. Write
  $a_w=(a_1,a_2)$, and write $\omega$ for the unique fundamental
  character of $I_{F_w}$ of niveau one.

Then $a_w\in\Wconj(\rbar|_{G_{F_w}})$.
\end{prop}
\begin{proof}
 
  Suppose first that $\rbar|_{G_{F_w}}$ is irreducible. Then the
  proof of Lemma 5.5 of \cite{geesavitttotallyramified} goes through
  unchanged, and gives the required result. So we may suppose that
  $\rbar|_{G_{F_w}}$ is reducible. In this case the proof of Lemma 5.4 of
  \cite{geesavitttotallyramified} goes through unchanged, and shows
  that we have \[\rbar|_{G_{F_w}}\cong
  \begin{pmatrix}
    \chibar_1&*\\0&\chibar_2
  \end{pmatrix}\]where
  $(\chibar_1\chibar_2)|_{I_K}=\omega^{a_1+a_2+e}$, and either
  $\chibar_1|_{I_K}=\omega^{a_1+z}$ or
  $\chibar_1|_{I_K}=\omega^{a_2+e-z}$ for some $1\le z\le e$, so we
  are in the situation of Section \ref{subsec: H^1_f}. Consider the
  extension class in $H^1(G_{F_w},\chibar_1\chibar_2^{-1})$
  corresponding to $\rbar|_{G_{F_w}}$. By Proposition \ref{prop:
    modular of some weight implies potentially BT lifts exist}, either
  $a_1-a_2=p-1$ and $\chibar_1\chibar_2^{-1}=\epsilonbar$, or this extension class is in $\Lflat$. In either case,
  by Theorem \ref{thm: crystalline equals flat} and Proposition
  \ref{prop: Lcrys in the exceptional case}, the extension class is in
  $\Lcrys$, so that  $a_w\in\Wconj(\rbar|_{G_{F_w}})$, as required.
\end{proof}
Combining this with Theorem 5.1.3 of \cite{blggU2}, we obtain our
final result.
\begin{thm}
  \label{thm: the main result, modular if and only if predicted}Let
  $F$ be an imaginary CM field with maximal totally real subfield
  $F^+$, and suppose that $F/F^+$ is unramified at all finite places,
  that every place of $F^+$ dividing $p$ splits completely in $F$,
  that $\zeta_p\notin F$, and that $[F^+:\Q]$ is even. Suppose that
  $p>2$, and that $\rbar:G_F\to\GL_2(\Fpbar)$ is an irreducible
  modular representation with split ramification such that
  $\rbar(G_{F(\zeta_p)})$ is adequate. Assume that for each place $w|p$
  of $F$, $F_w/\Qp$ is totally ramified.

 Let $a\in(\Z^2_+)_0^S$ be a Serre weight. Then
 $a_w\in\Wconj(\rbar|_{G_{F_w}})$ for all $w$ if and only if $\rbar$ is modular of
 weight $a$.
\end{thm}

\bibliographystyle{amsalpha}
\bibliography{geeliusavitt}

\end{document}